\documentclass[12pt]{article}
\usepackage{amsmath,amssymb,amsthm,amscd,graphicx}
\usepackage{latexsym}
\usepackage{enumerate}
\usepackage{amsbsy, amsfonts, textcomp}
\usepackage{makeidx}
\newtheorem{theorem}{Theorem}[]
\newtheorem{proposition}[theorem]{Proposition}
\newtheorem{lemma}[theorem]{Lemma}
\newtheorem{corollary}[theorem]{Corollary}

\theoremstyle{definition}

\newtheorem{example}[theorem]{Example}

\theoremstyle{remark}
\newtheorem{remark}[theorem]{Remark}
\newtheorem{question}[theorem]{Question}

\newcommand{\C}{\mathbb{C}}

\newcommand{\N}{\mathbb{N}}
\newcommand{\A}{\mathbb{A}}

\def \GL{\operatorname{GL}}

\def \max{\operatorname{max}}

\begin{document}

\begin{center}

\Large

{\bf On some nice polynomial automorphisms}

\normalsize

\vspace{1cm}

EL\.ZBIETA ADAMUS \\ Faculty of Applied Mathematics, \\ AGH University of Science and Technology \\
al. Mickiewicza 30, 30-059 Krak\'ow, Poland \\
e-mail: esowa@agh.edu.pl \\

\vspace{0.5cm}
PAWE\L \ BOGDAN \\ Faculty of Mathematics and Computer Science, \\ Jagiellonian University \\
ul. \L ojasiewicza 6, 30-348 Krak\'ow, Poland \\
e-mail: pawel.bogdan@uj.edu.pl \\

\vspace{0.5cm}
TERESA CRESPO \\ Departament d'\`{A}lgebra i Geometria, \\ Universitat de Barcelona \\
Gran Via de les Corts Catalanes 585, 08007 Barcelona, Spain \\
e-mail: teresa.crespo@ub.edu

\vspace{0.5cm}
ZBIGNIEW HAJTO \\ Faculty of Mathematics and Computer Science, \\ Jagiellonian University \\
ul. \L ojasiewicza 6, 30-348 Krak\'ow, Poland \\
e-mail: zbigniew.hajto@uj.edu.pl

\end{center}

\vspace{0.2cm}

\begin{abstract}
Given a polynomial endomorphism $F$ of the $n$-dimensional affine space $\A_K^n$ over a field $K$, we define a sequence of polynomial endomorphisms $(P_k)_{k \in \N}$ of $\A_K^n$ associated to $F$. We call $F$ nice if there exists an integer $m$ such that $P_m=0$ for some $m$. Then $F$ is invertible and $F^{-1}$ may be computed in terms of the endomorphisms $P_k, k=0,\dots ,m-1$. In this paper we study the class of nice polynomial automorphisms and obtain that this class is large and includes triangulable automorphisms and all linear cubic homogeneous polynomial automorphisms of nilpotence index up to 3. We prove as well that the nicety property is invariant under linear conjugation and determine that seven of the eight forms in Hubbers' classification of cubic homogeneous automorphisms in dimension 4 are nice and moreover the eighth one is a composition of nice maps.

\end{abstract}

\section{Introduction}

The Jacobian Conjecture originated in the question raised by Keller in \cite{K} on the invertibility of polynomial maps with Jacobian determinant equal to~1. The question is still open in spite of the efforts of many mathematicians. We recall in the sequel the precise statement of the Jacobian Conjecture and other results we shall use. We refer to \cite{E} for a detailed account of the research on the Jacobian Conjecture and related topics.

Let $K$ be a field and $K[X]=K[X_1,\dots,X_n]$ the polynomial ring in the variables $X_1,\dots,X_n$ over $K$. A \emph{polynomial map} is a map $F=(F_1,\dots,F_n):K^n \rightarrow K^n$ of the form

$$(X_1,\dots,X_n)\mapsto (F_1(X_1,\dots,X_n),\dots,F_n(X_1,\dots,X_n)),$$

\noindent where $F_i \in K[X], 1 \leq i \leq n$. The polynomial map $F$ is \emph{invertible} if there exists a polynomial map $G=(G_1,\dots,G_n):K^n \rightarrow K^n$ such that $X_i=G_i(F_1,\dots,F_n),  1 \leq i \leq n$. We shall call $F$ a \emph{Keller map} if the Jacobian matrix

$$J=\left(\dfrac{\partial F_i}{\partial X_j}\right)_{\substack{1\leq i \leq n \\ 1\leq j \leq n}}$$

\noindent has determinant equal to a nonzero constant. Clearly an invertible polynomial map $F$ is a Keller map.

\vspace{0.5cm}
\noindent {\bf Jacobian Conjecture.} {\it Let $K$ be a field of characteristic zero. A Keller map $F:K^n \rightarrow K^n$ is invertible.}

\vspace{0.5cm}
In the sequel, $K$ will always denote a field of characteristic $0$.
For $F=(F_1,\dots,F_n) \in K[X]^n$, we define the \emph{degree} of $F$ as $\deg F= \max \{\deg F_i : 1\leq i \leq n\}$. It is known that if $F$ is a polynomial automorphism, then $\deg F^{-1} \leq (\deg F)^{n-1}$ (see \cite{BCW} or \cite{RW}).

The Jacobian conjecture for quadratic maps was proved by Wang in \cite{W}.  We state now the reduction of the Jacobian conjecture to the case of maps of third degree (see \cite{BCW}, \cite{Y}, \cite{D1} and \cite{D2}).

\begin{proposition}\label{red} \begin{enumerate}[a)] \item (Bass-Connell-Wright-Yagzhev) Given a Keller map $F:K^n \rightarrow K^n$, there exists a Keller map $\widetilde{F}:K^N \rightarrow K^N$, $N\geq n$ of the form $\widetilde{F}=Id+H$, where $H(X)$ is a cubic homogeneous map and having the following property: if $\widetilde{F}$ is invertible, then $F$ is invertible too.
\item (Dru\.zkowski) The cubic part $H$ may be chosen of the form

$$\left( (\sum_{j=1}^N a_{1j} X_j)^3, \dots,(\sum_{j=1}^N a_{Nj} X_j)^3 \right)$$

\noindent and with the matrix $A=(a_{ij})_{\substack{1\leq i \leq N \\ 1\leq j \leq N}}$ satisfying $A^2=0$.

\end{enumerate}

\end{proposition}

We recall now the sufficient condition for the invertibility of a polynomial map we gave in \cite{ABCH} Corollary 5.

\begin{proposition}\label{suff} Given a polynomial map
$F:K^n \rightarrow K^n$, we define for $i=1,\dots,n$ the following
sequence of polynomials in $K[X]$,

$$\begin{array}{lll} P_0^i(X_1,\dots,X_n) &= & X_i, \\
P_1^i(X_1,\dots,X_n) &=& F_i-X_i,
\end{array}$$

\noindent and, assuming $P_{k-1}^i$ is defined,

$$P_k^i(X_1,\dots,X_n) = P_{k-1}^i(F_1,\dots,F_n)-P_{k-1}^i(X_1,\dots,X_n).$$

If there
exists an integer $m$ such that $P_{m}^i = 0$, for all $i=1,\dots,
n$, then $F$ is invertible and the inverse $G$ of $F$ is given by

$$G_i(Y_1,Y_2,\dots,Y_n)= \sum_{l=0}^{m-1} (-1)^l P_l^i(Y_1,Y_2,\dots,Y_n), \, 1 \leq i \leq n.$$

\end{proposition}

We shall say that the polynomial map $F$ is \emph{nice} if there
exists an integer $m$ such that $P_{m}^i = 0$, for all $i=1,\dots,
n$. The aim of this paper is to prove several properties of nice polynomial automorphisms, namely that nicety is invariant under linear conjugation, that all triangulable polynomial automorphisms are nice and that the inverse of a nice polynomial automorphism is nice. We also prove that if $f$ and $F$ form a Gorni-Zampieri pairing and $F$ is nice, then $f$ is nice and that all linear cubic polynomial maps of nilpotence index $\leq 3$ are nice.

\section{Triangulable polynomial automorphisms}
\begin{theorem}\label{conj}
Let $F:K^n \rightarrow K^n$ be a polynomial map of the form

$$\left\{ \begin{array}{lll} F_1(X_1,\dots,X_n)&=& X_1+H_1(X_1,\dots,X_n) \\ F_2(X_1,\dots,X_n)&=& X_2+H_2(X_1,\dots,X_n)\\ & \vdots & \\  F_n(X_1,\dots,X_n)&=& X_n+H_n(X_1,\dots,X_n) , \end{array} \right.$$

\noindent where $H_i(X_1,\dots,X_n)$ is a polynomial in
$X_1,\dots,X_n$ of lower degree $\geq 2$, $1\leq i \leq n$. Let $T
\in \GL(n,K)$ and $\widetilde{F}:= T^{-1} \circ F \circ T$. If $F$
is nice, then $\widetilde{F}$ is nice.
\end{theorem}

\noindent {\it Proof.} We have clearly $\widetilde{F}:= T^{-1}
\circ F \circ T= T^{-1} \circ (I+H) \circ T=I+\widetilde{H}$,
where $\widetilde{H}=T^{-1} \circ H \circ T$ has lower degree
$\geq 2$. Let us consider the polynomial sequences $P_k^i$ and
$\widetilde{P}_k^i$ defined by

$$\begin{array}{lll} P_0^i(X_1,\dots,X_n) &= & X_i, \\
P_1^i(X_1,\dots,X_n) &=& F_i-X_i, \\
& \vdots & \\ P_k^i(X_1,\dots,X_n) &=&
P_{k-1}^i(F_1,\dots,F_n)-P_{k-1}^i(X_1,\dots,X_n)\\
& \vdots & \end{array}$$

\noindent and

$$\begin{array}{lll} \widetilde{P}_0^i(X_1,\dots,X_n) &= & X_i, \\
\widetilde{P}_1^i(X_1,\dots,X_n) &=& \widetilde{F}_i-X_i, \\
& \vdots & \\ \widetilde{P}_k^i(X_1,\dots,X_n) &=&
\widetilde{P}_{k-1}^i(\widetilde{F}_1,\dots,\widetilde{F}_n)-\widetilde{P}_{k-1}^i(X_1,\dots,X_n)\\
& \vdots &
\end{array}$$

\noindent and the polynomial maps $P_k:=(P_k^1,\dots,P_k^n)$ and
$\widetilde{P}_k:=(\widetilde{P}_k^1,\dots,\widetilde{P}_k^n)$. We
shall prove by induction on $k$

$$\widetilde{P}_k=T^{-1} \circ P_k \circ T, \text{\ for all integer\ } k\geq 1.$$

\noindent Then, if $P_m=0$ for some integer $m$, also
$\widetilde{P}_m=0$. For $k=1$, we have
$\widetilde{P}_1=\widetilde{H}=T^{-1} \circ H \circ T=T^{-1} \circ
P_1 \circ T.$ Let us assume $\widetilde{P}_k=T^{-1} \circ P_k
\circ T$. We have then

$$\begin{array}{lll} \widetilde{P}_{k+1}&=&  \widetilde{P}_k \circ \widetilde{F}
-\widetilde{P}_k \\ &=& (T^{-1} \circ P_k \circ T) \circ (T^{-1}
\circ F \circ T)-T^{-1} \circ P_k \circ T \\ &=& T^{-1} \circ P_k
\circ F \circ T-T^{-1} \circ P_k \circ T \\ &=& T^{-1} \circ (P_k
\circ F - P_k) \circ T \\ &=& T^{-1} \circ P_{k+1} \circ
T.\end{array}$$

\hfill $\Box$

\begin{theorem}\label{tri} Let $F$ be a triangular polynomial map in dimension $n$, i.e. a polynomial map of the form

$$\left\{ \begin{array}{lll} F_1(X_1,\dots,X_n)&=& X_1+H_1(X_2,\dots,X_n) \\ F_2(X_1,\dots,X_n)&=& X_2+H_2(X_3,\dots,X_n)\\ & \vdots & \\  F_{n-1}(X_1,\dots,X_n)&=& X_{n-1}+H_{n-1}(X_n) \\  F_n(X_1,\dots,X_n)&=& X_n. \end{array} \right.$$

Then $F$ is nice.
\end{theorem}

\noindent {\it Proof.} Let $H \in K[X_2,\dots,X_n]$ have degree $d$ in $X_2$ and define

$$\begin{array}{lll} P_1(X_2,\dots,X_n) &= & H, \\
P_k(X_2,\dots,X_n) &=& P_{k-1}(F_2,\dots,F_n)-P_{k-1}(X_2,\dots,X_n).
\end{array}$$

\noindent Since the polynomials $H_j$, for $j>1$, have degree $0$ in $X_2$, we obtain that the degree in $X_2$ of $P_k$ is $d-k+1$, hence the degree in $X_2$ of $P_{d+1}$ is $0$. Taking $H=H_1$, we obtain that $P_{d+1}^1$ is a polynomial in $X_3,\dots,X_n$. If $P_{d+1}^1$ has degree $d'$ in $X_3$, taking $H=P_{d+1}^1$ and the polynomial map $(F_2,\dots, F_n)$, we obtain that $P_{d+d'+1}^1$ is a polynomial in $X_4,\dots,X_n$. By iteration, we obtain $P_m^1=0$, for some $m$, hence $F$ is nice. \hfill $\Box$

\vspace{0.5cm}

Theorems \ref{conj} and \ref{tri} together give the following corollary.

\begin{corollary}
A triangulable polynomial automorphism is nice.
\end{corollary}

\begin{example} We consider the Nagata automorphism defined by

$$\left\{ \begin{array}{ccl} F_1 &=& X_1-2pX_2-p^2 X_3 \\ F_2 &=& X_2+pX_3 \\ F_3 &=& X_3 \end{array} \right. $$

\noindent where $p=X_1 X_3+X_2^2$. It is known that the Nagata automorphism is not tame. However it is nice. Indeed, by calculation, we obtain $P_3^1=0, P_2^2=0$. \end{example}

\section{Products of nice automorphisms}

\begin{theorem} If $F$ is a nice polynomial automorphism of $K^n$, then $F^{-1}$ is nice.
\end{theorem}

\noindent {\it Proof.} Let $F(X)$ be nice and $G(Y)$ be its inverse. Let the polynomials $P_k^i$ be defined by

$$\begin{array}{lll} P_0^i(X_1,\dots,X_n) &= & X_i, \\
P_1^i(X_1,\dots,X_n) &= & F_i-X_i, \\
P_k^i(X_1,\dots,X_n) &=& P_{k-1}^i(F_1,\dots,F_n)-P_{k-1}^i(X_1,\dots,X_n)
\end{array}$$

\noindent and assume $P_m^i=0$. Then by Proposition \ref{suff}, we have

$$X_i= \sum_{l=0}^{m-1} (-1)^l P_l^i(F)$$

\noindent hence

\begin{equation}\label{inv} G_i(Y)= \sum_{l=0}^{m-1} (-1)^l P_l^i(Y).
\end{equation}

Let now the polynomials $Q_k^i$ be defined by

$$\begin{array}{lll} Q_0^i(Y_1,\dots,Y_n) &= & Y_i, \\
Q_1^i(Y_1,\dots,Y_n) &= & G_i-Y_i, \\
Q_k^i(Y_1,\dots,Y_n) &=& Q_{k-1}^i(G_1,\dots,G_n)-Q_{k-1}^i(Y_1,\dots,Y_n).
\end{array}$$

We shall compute the polynomials $Q_k^i$ using (\ref{inv}). We have

$$\begin{array}{l} Q_1^i(Y) =  G_i-Y_i = \sum_{l=1}^{m-1} (-1)^l P_l^i(Y) \\ Q_2^i(Y) =   Q_1^i(G)- Q_1^i(Y) = \sum_{l=1}^{m-1} (-1)^l ( P_l^i(G)- P_l^i(Y)) \end{array}$$

\noindent Now $P_k^i(X) = P_{k-1}^i(F)-P_{k-1}^i(X)$ implies $P_k^i(G) = P_{k-1}^i(Y)-P_{k-1}^i(G)$, therefore

$$ Q_2^i(Y) = \sum_{l=2}^{m-1} (-1)^l  P_l^i(G)$$

\noindent and, by induction, we obtain

$$ Q_k^i(Y) = \sum_{l=k}^{m-1} (-1)^l  P_l^i(G^{k-1}),$$

\noindent where $G^{k-1}$ denotes composition of $G$ with itself $k-1$ times. We conclude then that $Q_m^i(Y)=0$.
\hfill $\Box$

\begin{remark} The set of nice polynomial automorphisms is not a subgroup of the group of polynomial automorphisms as can be shown by considering the following example. Let $F_1$ and $F_2$ be the polynomial automorphisms of $K^2$ defined by $X \mapsto (X_1+X_2^3,X_2)$ and $X \mapsto (X_1,X_2+X_1^2)$ respectively. It is easy to check that $F_1$ and $F_2$ are nice. However $F:=F_1 \circ F_2$ is not nice. We have $F(X)=(X_1+(X_1^2+X_2)^3,X_2+X_1^2)$. Let $Q_k$ denote the term of lower degree of $P_k^1$. Taking into account that $Q_{k+1}$ is the term of lower degree of

$$\dfrac{\partial Q_k}{\partial X_1} (X_1^2+X_2)^3 + \dfrac{\partial Q_k}{\partial X_2} X_1^2$$

\noindent and that, if $Q$ is a homogeneous polynomial, the term of lower degree of \linebreak $\dfrac{\partial Q}{\partial X_1} (X_1^2+X_2)^3 + \dfrac{\partial Q}{\partial X_2} X_1^2$ is

$$\begin{array}{lll} \dfrac{\partial Q}{\partial X_2} X_1^2 & \text{\ if \ } & \dfrac{\partial Q}{\partial X_2} \neq 0 \\[10pt]
\dfrac{\partial Q}{\partial X_1} X_2^3 & \text{\ if \ } & \dfrac{\partial Q}{\partial X_2} = 0 \end{array} $$

\noindent we obtain the following recursive formula for $Q_k$, which shows $P_k^1 \neq 0$ for all $k$. Let us write $k= 4m+r$, with $r=1,2,3,4$ and define $\mu(r)$ by $\mu(1)=1, \mu(2)=3, \mu(3)=\mu(4)=6$. Then

$$Q_k= \mu(r) 6^m \prod_{j=0}^m (6+5(j-1)) X_1^{2r-2+5m} X_2^{4-r}.$$

\end{remark}

\begin{remark}[{\it Hubbers' classification}]
Hubbers \cite{Hu} classified up to conjugation all cubic homogeneous polynomial maps in dimension 4. He obtained eight conjugation classes. We have checked that all but the last are nice automorphisms. Let us now consider the eighth class in Hubbers classification. It is represented by

$$\left\{\begin{array}{ccl} F_1 & = & X_ 1 \\ F_2 & = & X_ 2 - X_1^3/3 \\ F_3 & = & X_3-X_1^2X_2-e_3X_1X_2^2+g_4X_1X_2X_3-k_3X_2^3+m_4X_2^2X_3+g_4^2X_2^2X_4 \\
F_4 &=& X_4-X_1^2X_3-e_4X_1X_2^2-2m_4X_1X_2X_3/g_4-g_4X_1X_2X_4-k_4X_2^3\\ && -m_4^2X_2^2X_3/g_4^2-m_4X_2^2X_4 \end{array}.\right.$$

\noindent We have $F=G\circ H$, where $H$ and $G$ defined by

$$\left\{\begin{array}{ccl} H_1 & = & X_ 1 \\ H_2 & = & X_ 2 \\ H_3 & = & F_3 \\
H_4 &=& F_4 \end{array}\right. \text{\ \ and \ \ } \left\{\begin{array}{ccl} G_1 & = & X_ 1 \\ G_2 & = & F_2 \\ G_3 & = & X_3 \\ G_4 & = & X_4 \end{array}\right.$$

\noindent are both nice.
\end{remark}

\begin{remark} It is known that all cubic homogeneous polynomial automorphisms in dimension $\leq 3$ are triangulable (see \cite{E}, Proposition 7.1.1), hence they are nice. We obtain then that the subset of nice cubic homogeneous polynomial automorphisms in dimension $\leq 4$ generate the whole set of cubic homogeneous polynomial automorphisms.
\end{remark}

\begin{question} The above remark suggest the following question: Is every polynomial automorphism a product of nice ones?
\end{question}

\section{Gorni-Zampieri pairing}

We shall now prove that if $f$ and $F$ form a Gorni-Zampieri pair (see \cite{E} \S 6.4) and $F$ is nice then $f$ is nice.

Let $f= K^n\rightarrow K^n$ be a cubic homogeneous map and $F: K^N \rightarrow K^N$ a cubic linear map, with $N>n$.
Let $x =(x_1,\ldots, x_n) \in K^n$ and $X=(X_1, \ldots, X_N) \in K^N$.

For $f=(f_1,\ldots, f_n)$, where $f_i \in K[x]$ we define $p_t(x)=(p_t^1(x), \ldots, p_t^n(x))^T$ in the following way

$$\begin{array}{lll} p_0^i(x)&=& x_i \\ p_1^i(x)&=&f_i(x)-x_i \\ p_2^i(x)&=&p_1^i(f(x))-p_1^i(x) \\ & \vdots & \\ p_{t+1}^i(x)&=& p_t^i(f(x))-p_t^i(x) \end{array} $$

For $F=(F_1,\ldots, F_N)$, where $F_i \in K[X]$ we define $P_t(X)=(P_t^1(X), \ldots, P_t^n(X))^T$ in the following way

\[ \begin{array}{lll} P_0^i(X)&=& X_i\\ P_1^i(X)&=&F_i(X)-X_i\\ P_2^i(X)&=&P_1^i(F(x))-P_1^i(X)\\ &\vdots &  \\ P_{t+1}^i(X)&=& P_t^i(F(X))-P_t^i(X)\end{array} \]

\begin{lemma}\label{pair}
Let $f$ and $F$ be as above.
Assume that $f$ and $F$ are paired through the matrices $B \in M_{n \times N}(K)$ and $C \in M_{N \times n}(K)$.
Then $p_t$ and $P_t$ are paired through the same matrices $B$ and $C$, for every $t\geq 0$.
\end{lemma}

\noindent \textit{Proof.} Our main assumption is that $F$ and $f$ are paired, so in particular $f(x)=BF(X)$, where $X=Cx$.
Note that since $BC=I_n$, then $X=Cx$ implies $x=BX$.

The proof goes by induction on $t$. For $t=0$ we have $p_0(x)=x=BX$ and $P_0(X)=X$, so $p_0(x)=BP_0(X)$.
Similarly for $t=1$ we have
\[p_1(x)=p_0(f(x))-p_0(x)=f(x)-x,\]
\[P_1(X)=P_0(F(X))-P_0(X)=F(X)-X.\]
So
\[p_1(x)=f(x)-x=BF(X)-BX=B(F(X)-X)=BP_1(X).\]

Let us assume that for a given $t$ we have
\begin{equation}
 p_t(x)=BP_t(X)
 \label{one}
\end{equation}
or equivalently
\begin{equation}
 p_t(BX)=BP_t(X) \qquad \forall X \in K^N.
 \label{two}
\end{equation}
We will prove that
\[
 p_{t+1}(x)-BP_{t+1}(X)=0.
\]
We have that
\[  p_{t+1}(x)-BP_{t+1}(X)=p_t(f(x)) - p_t(x) - [BP_t(F(X))-BP_t(X)]=^{(\ref{one})} \]
\[ =p_t(f(x)) -BP_t(F(X))=p_t(BF(X))-BP_t(F(X))=^{(\ref{two})}0\]
\hfill $\Box$

\begin{corollary}\label{GZ}
 Let $f$ and $F$ be as above. If $f$ and $F$ are paired, then if $F$ is nice, $f$ is also nice.
\end{corollary}

\noindent {\it Proof.} If $f$ and $F$ are paired through the matrices $B$ and $C$, Lemma \ref{pair} gives  $p_{t}(x)=BP_{t}(X)$, for all $t\geq 0$. Then if $P_m=0$, for some integer $m$, we have $p_m=0$.
\hfill $\Box$

\begin{example} In \cite{ABCH} Example 7, we proved that the polynomial automorphism $f$ of $K^4$ defined by

$$\left\{ \begin{array}{lll} f_1 &=& x_1+px_4 \\ f_2 &=& x_2-px_3 \\f_3 &=& x_3+x_4^3 \\ f_4 &=& x_4 \end{array}, \right.$$

\noindent where $p:=x_1x_3+x_2x_4$, is not nice. We consider the polynomial automorphism $F$ of $K^{16}$ defined by

$$\left\{ \begin{array}{rll} F_1 &=& X_1+(L_1+X_{15}+X_{16})^3 \\ F_2 &=& X_2+(L_1-X_{15}-X_{16})^3 \\ F_3 &=& X_3+(L_1+X_{15}-X_{16})^3 \\ F_4 &=& X_4+(L_1-X_{15}+X_{16})^3 \\
F_5 &=& X_5+(L_2+X_{16})^3 \\ F_6 &=& X_6+(L_2-X_{16})^3 \\ F_7 &=& X_7+L_2^3 \\ F_8 &=& X_8+(L_1+X_{15})^3 \\ F_9 &=& X_9+(L_1-X_{15})^3 \\ F_{10} &=& X_{10}+L_1^3 \\
F_{11} &=& X_{11}+(L_2+X_{15}+X_{16})^3 \\ F_{12} &=& X_{12}+(L_2-X_{15}-X_{16})^3 \\ F_{13} &=& X_{13}+(L_2+X_{15}-X_{16})^3 \\ F_{14} &=& X_{14}+(L_2-X_{15}+X_{16})^3
\\ F_{15} &=& X_{15}+X_{16}^3\\ F_{16} &=& X_{16} \end{array}\right.
$$

\noindent where

$$\begin{array}{lll} L_1 &=& (X_1+X_2-X_3-X_4+4X_5+4X_6-8X_7)/24 \\ L_2 &=& (-4X_8-4X_9+8X_{10}-X_{11}-X_{12}+X_{13}+X_{14})/24 \end{array}$$

\noindent which is paired with $f$ through the matrices

$$B=\dfrac 1 {24}  \left( \begin{array}{rrrrrrrrrrrrrrrr} 1&1&-1&-1& 4&4&-8&0&0&0&0&0&0&0&0&0 \\ 0&0&0&0&0&0&0&-4&-4&8&-1&-1&1&1&0&0\\ 0&0&0&0&0&0&0&0&0&0&0&0&0&0&24&0 \\ 0&0&0&0&0&0&0&0&0&0&0&0&0&0&0&24 \end{array} \right)$$

\noindent and

 $$C= \left( \begin{array}{rrrrrrrrrrrrrrrr} 0&0&0&0&0&0&-3&0&0&0&0&0&0&0&0&0 \\0&0&0&0&0&0&0&-6&0&0&0&0&0&0&0&0 \\ 0&0&0&0&0&0&0&0&0&0&0&0&0&0&1&0 \\0&0&0&0 &0&0&0&0&0&0&0&0&0&0&0&1 \end{array} \right)^T.$$

\noindent Applying Corollary \ref{GZ}, we obtain that $F$ is not nice. Let us note that the Jacobian matrix of $F-Id$ has nilpotence index equal to 5.

\end{example}

\section{Linear cubic polynomial maps}

We shall see now that all polynomial maps of the form $F=X+H$, where each component of $H$ is the cube of a linear form and the Jacobian matrix of $H$ has nilpotence index $\leq 3$, are nice.

\begin{theorem}
Let $F:K^n \rightarrow K^n$ be a polynomial map of the form

$$\left\{ \begin{array}{lll} F_1(X_1,\dots,X_n)&=& X_1+H_1(X_1,\dots,X_n) \\ F_2(X_1,\dots,X_n)&=& X_2+H_2(X_1,\dots,X_n)\\ & \vdots & \\  F_n(X_1,\dots,X_n)&=& X_n+H_n(X_1,\dots,X_n) , \end{array} \right.$$

\noindent where $H_i(X_1,\dots,X_n)=L_i(X_1,\dots,X_n)^3$ and $L_i(X_1,\dots,X_n)=a_{i1}X_1+\dots +a_{in} X_n$, $1\leq i \leq n$. Let $JH$ denote the jacobian matrix  of $H:=(H_1,\dots,H_n)$.

\begin{enumerate}[1)]
\item
If $(JH)^2=0$, then $F$ is nice. Moreover $F^{-1}=X-H$.
\item
If $(JH)^3=0$, then $F$ is nice. Moreover $F^{-1}$ has degree at most 9.
\end{enumerate}
\end{theorem}

\noindent {\it Proof.} 1) Since $H$ is homogeneous, $(JH)^2=0$ implies $(JH)H=0$. We have then $\sum_{j=1}^n \dfrac {\partial H_i}{\partial X_{j}}H_j=0,  1\leq i \leq n$. We shall prove

\begin{equation}\label{eq} \sum_{1\leq j_1,\dots,j_k \leq n} \dfrac {\partial^k H_i}{\partial X_{j_1} \dots \partial X_{j_k}}H_{j_1} \dots H_{j_k}=0, 1\leq i \leq n,
\end{equation}

\noindent for all $k\geq 1$, by induction on $k$. For $k=1$, it is true by hypothesis. Let us assume that it is true for $k$. By applying $\partial /\partial X_j$ to (\ref{eq}), we obtain

$$\sum_{1\leq j_1,\dots,j_k \leq n} (\dfrac {\partial^{k+1} H_i}{\partial X_{j_1} \dots \partial X_{j_k}\partial X_j}H_{j_1} \dots H_{j_k}+\sum_{l=1}^k \dfrac {\partial^k H_i}{\partial X_{j_1} \dots \partial X_{j_k}}H_{j_1}\dots \dfrac {\partial H_{j_l}}{\partial X_j} \dots H_{j_k})=0.$$

\noindent Multiplying by $H_j$, summing from $j=1$ to $n$ and taking into account the case $k=1$ with $i=j_1,\dots,j_k$, we obtain (\ref{eq}) for $k+1$. Lemma 6 in \cite{ABCH} gives then $P_2^i=0, 1\leq i \leq n$, hence $F$ is nice and $F^{-1}=X-H$.

2) We consider, for $i=1,\dots,n$,

$$\begin{array}{ccl} P_0^i(X) &:=& X_i \\ P_1^i(X) &:=& F_i-X_i=H_i \\ P_2^i(X) &:=& H_i(F)-H_i(X) \\ P_3^i(X) &:=& P_2^i(F)-P_2^i(X) \\ P_4^i(X) &:=& P_3^i(F)-P_3^i(X) \\P_5^i(X) &:=& P_4^i(F)-P_4^i(X). \end{array} $$

\noindent Our aim is to prove $P_5^i=0, \, \forall i=1,\dots,n,$ and $\deg (P_k^i) \leq 9, \, \forall i=1,\dots,n, \, k=0,\dots, 4$. Applying \cite{ABCH} Corollary 5, this implies that the map $F$ is invertible and $\deg F^{-1} \leq 9$.

Since  $(JH)^3=0$, we have

\[
\sum_{j,k=1}^n \dfrac{\partial H_i}{\partial X_j}\dfrac{\partial H_j}{\partial X_k}\dfrac{\partial H_k}{\partial X_l}=0, \quad \forall \, i, l =1,\dots,n.
\]

\noindent i.e.

\[
\sum_{j,k=1}^n (3L_i^2 a_{ij}) (3L_j^2 a_{jk}) (3L_k^2 a_{kl}) =27 L_i^2 (\sum_{j,k=1}^n  a_{ij}L_j^2 a_{jk} L_k^2 a_{kl})=0, \quad \forall \, i, l =1,\dots,n.
\]

\noindent which implies

\begin{equation}\label{eq3}
\sum_{j,k=1}^n a_{ij}  a_{jk} a_{kl} L_j^2 L_k^2  =0, \quad \forall \, i, l =1,\dots,n,
\end{equation}

\noindent since $L_i=0$ implies $a_{ij}=0$, for all $j$, hence (\ref{eq3}). Applying $\dfrac{\partial}{\partial X_m}$, we obtain

\begin{equation}\label{eq4}
\sum_{j,k=1}^n a_{ij} a_{jk} a_{kl}(a_{jm} L_jL_k^2+ a_{km}L_j^2L_k ) =0, \quad \forall \, i, l, m =1,\dots,n.
\end{equation}

\noindent Further, from (\ref{eq3}), we obtain

\begin{equation}\label{eq8}
\sum_{j,k=1}^n a_{ij} a_{jk}L_j^2L_k^3 =0, \quad \forall \, i=1,\dots,n
\end{equation}

\noindent and, from (\ref{eq4}),

$$\sum_{j,k=1}^n a_{ij} a_{jk} (a_{jm}L_jL_k^3+a_{km}L_j^2L_k^2 ) =0, \quad \forall \, i, m =1,\dots,n.$$

\noindent Using (\ref{eq3}) again, we have

\begin{equation}\label{eq9}
\sum_{j,k=1}^n a_{ij} a_{jk} a_{jm} L_jL_k^3 =0, \quad \forall \, i, m =1,\dots,n.
\end{equation}

\noindent Applying $\partial/\partial X_p$ to (\ref{eq9}), we obtain

$$\sum_{j,k=1}^n a_{ij} a_{jk} a_{jm} (a_{jp} L_k^3+3a_{kp} L_j L_k^2) =0$$

\noindent i.e.

\begin{equation}\label{eq10}
\sum_{j,k=1}^n a_{ij} a_{jk} a_{jm} a_{jp} L_k^3= - 3\sum_{j,k=1}^n a_{ij} a_{jk} a_{jm} a_{kp} L_j L_k^2
\end{equation}

Computing $P_2^i, P_3^i, P_4^i$ and $P_5^i$ with the formulas given in the proof of Lemma 6 in \cite{ABCH} and taking into account equalities (\ref{eq8}), (\ref{eq9}) and (\ref{eq10}), we obtain

$$\begin{array}{lll} P_2^i &=&  3L_i^2 \sum_{j=1}^n a_{ij} L_j^3+3L_i \sum_{1\leq j,k \leq n} a_{ij} a_{ik} L_j^3 L_k^3+ \sum_{1\leq j,k,l \leq n}  a_{ij} a_{ik} a_{il} L_j^3 L_k^3 L_l^3 \\ P_3^i &=& 6L_i \sum_{j,k=1}^n a_{ij}a_{ik} L_j^3 L_k^3+\sum_{j,k,l=1}^n 9 a_{ij} a_{ik} a_{il} L_j^3L_k^3L_l^3 \\ P_4^i &=& \sum_{j,k,l=1}^n 6 a_{ij}a_{ik}a_{il}  L_j^3 L_k^3 L_l^3 \\ P_5^i &=& 0.
\end{array} $$
\hfill $\Box$

 \vspace{0.5cm} \noindent {\bf Acknowledgments.} E.
Adamus acknowledges support of the Polish Ministry of Science and
Higher Education. T. Crespo and Z. Hajto acknowledge support of
grant MTM2012-33830, Spanish Science Ministry.

\end{document}